\documentclass[12pt]{article}

\begin{document}
\bibliographystyle{plain}

\thispagestyle{empty}
\setcounter{page}{0}

{\large Guszt\'av MORVAI and  Benjamin WEISS}
 
\vspace {2cm}

{\Large Prediction for Discrete  Time Series}

\vspace {2cm}

{\large Probab. Theory Related Fields  132  (2005),  no. 1, 1--12.}

\vspace {2cm}

\begin{abstract}
Let $\{X_n\}$ be a stationary and ergodic time series taking values from a finite or countably infinite 
set ${\cal X}$. 
Assume that the distribution of the process is otherwise unknown. 
We propose a sequence of stopping times $\lambda_n$ 
along which  we will be able to estimate the conditional probability 
$P(X_{\lambda_n+1}=x|X_0,\dots,X_{\lambda_n} )$
from data segment $(X_0,\dots,X_{\lambda_n})$
in a pointwise consistent way for a restricted class of stationary and ergodic  finite or countably infinite alphabet 
time series which includes among others all stationary and ergodic 
finitarily Markovian  processes. 
If the stationary and ergodic process turns out to be   finitarily Markovian (among others, all stationary and ergodic Markov chains 
are included in this class)    
then $ \lim_{n\to \infty} {n\over \lambda_n}>0$ almost surely. 
If the stationary and ergodic process turns out to possess   finite entropy rate  
 then  $\lambda_n$ is upperbounded by a polynomial, eventually almost surely.  
\end{abstract}



\pagebreak

\section{Introduction}

\noindent
Bailey \cite{Bailey76} and Ryabko \cite{Ryabko88} considered the problem of estimating the conditional probability 
$P(X_{n+1}=1|X_0,\dots,X_n)$ for  binary time series. They showed that one cannot estimate this 
quantity from the data $(X_0,\dots,X_n)$ such that the difference tends to zero 
almost surely 
as $n$ increases, for all 
 stationary and ergodic binary time series.

It is well known, that if one  knows in advance that the process is Markov with arbitrary (unknown)  order, 
then one can estimate the order 
(c.f. Csisz\'ar and Shields \cite{CsSh00}, Csisz\'ar \cite{Cs02}), 
and using this estimate for the order, one can count 
empirical averages of blocks with lengths one plus the order 
for estimating  $P(X_{n+1}=1|X_0,\dots,X_n)$ in a pointwise consistent way. 
In the present paper we will consider the case when it is
  not known in advance if the process is Markov or not.

Morvai \cite{Mo00} exhibited 
a sequence of stopping times  $\eta_n$ 
such that $P(X_{\eta_n+1}=1|X_0,\dots,X_{\eta_n} )$ can be 
estimated from data segment $(X_0,\dots,X_{\eta_n})$ in a pointwise consistent way, 
that is,  the error vanishes  as $n$ increases.  
The disadvantage of that scheme is that the stopping times grow very fast. 
Another, more reasonable scheme was proposed by  
Morvai and Weiss \cite{MW03} for a subclass  of stationary and ergodic 
binary time series. 
There  the stopping times still grow exponentially, though not so fast as in Morvai \cite{Mo00}.

Bailey \cite{Bailey76} proved that there is no test for the Markov property, 
that is, there is no algorithm which could tell you eventually if the process is 
Markov with any  order 
or not,  over all stationary and ergodic binary time series.

In this paper discrete (finite or countably infinite) alphabet  stationary and ergodic processes are treated. 
We   
propose a much denser (compared to Morvai and Weiss \cite{MW03})
 sequence of  stopping times $\lambda_n$  
 along which   we will be able to estimate 
$P(X_{\lambda_n+1}=x|X_0,\dots,X_{\lambda_n} )$ from samples $(X_0,\dots,X_{\lambda_n})$ in a pointwise consistent way
 for those processes whose conditional distribution is almost surely continuous
(see the precise definition below). This class includes all Markov processes with arbitrary order and the much wider class of finitarily Markovian processes. 
Despite  Bailey's result, for the proposed stopping times $\lambda_n$, 
if the stationary and ergodic process turns out to be finitarily Markovian 
(which includes all stationary and ergodic Markov chains with arbitrary order)  
then $ \lim_{n\to \infty} {n\over \lambda_n}>0$ almost surely. 
If the stationary and ergodic process  turns out to possess   finite entropy rate  then   
 $\lambda_n$ is upperbounded by a polynomial, eventually almost surely.

\section{The Proposed Algorithm}

\smallskip
\noindent
Let $\{X_n\}_{n=-\infty}^{\infty}$ be a stationary and ergodic time series taking values from a 
discrete (finite or countably infinite)  alphabet 
${\cal X}$. (Note that all stationary time series $\{X_n\}_{n=0}^{\infty}$ 
can be thought to be a 
two sided time series, that is, $\{X_n\}_{n=-\infty}^{\infty}$. )  
For notational convenience, let $X_m^n=(X_m,\dots,X_n)$,
where $m\le n$. Note that if $m>n$ then $X_m^n$ is the empty string.

\bigskip
\noindent
For $k\ge 1$, let $1\le l_k\le k$ be a nondecreasing unbounded sequence of integers, that is, 
$1=l_1\le l_2\dots$ and $\lim_{k\to\infty}l_k=\infty$.

\smallskip
\noindent
Define auxiliary  stopping times 
( similarly to 
Morvai and Weiss \cite{MW03}) 
as follows.  Set $\zeta_0=0$. 
For $n=1,2,\ldots$,  let
\begin{equation}\label{defzeta}
\zeta_n=\zeta_{n-1}+
\min\{t>0 : X_{\zeta_{n-1}-(l_n-1)+t}^{\zeta_{n-1}+t}=X_{\zeta_{n-1}-(l_n-1)}^{\zeta_{n-1}}\}.
\end{equation}

\smallskip
\noindent
Among other things, using  $\zeta_n$ and  $l_n$ we can  define a very useful process 
$\{ {\tilde X}_n\}_{n=-\infty}^{0}$  as a function of $X_0^{\infty}$ as follows. 
Let $J(n)=\min\{j\ge 1: \ l_{j+1}>n\}$ and  
define  
\begin{equation}
\label{defprocesses}
{\tilde X}_{-i}=X_{\zeta_{J(i)}-i} \ \ \mbox{for $i\ge 0$.}   
\end{equation}
As we will see in the proof of the Theorem, the $\{{\tilde X}\}_{n=-\infty}^0$ has the same 
distribution as the original process. 
For notational convenience let 
  $p_k(x_{-k}^0)$ and $p_k(y|x^0_{-k})$ denote the  distribution $P(X_{-k}^{0}=x_{-k}^{0})$ and  
the conditional distribution $P(X_1=y|X^0_{-k}=x^0_{-k})$, respectively. 

\bigskip
\noindent
{\bf Definition~$1$.} For a stationary  time series $\{X_n\}$ the (random) length $K(X^0_{-\infty})$ 
 of the memory of the sample path    
 $X^0_{-\infty}$ 
is the smallest possible $0\le K<\infty$ such that 
for all $i\ge 1$, all $y\in {\cal X}$, all $z^{-K}_{-K-i+1}\in {\cal X}^{i}$ 
$$
p_{K-1}(y|X^0_{-K+1})=p_{K+i-1}(y|z^{-K}_{-K-i+1},X^0_{-K+1})
$$ 
provided $p_{K+i}(z^{-K}_{-K-i+1},X^0_{-K+1},y)>0$,
and $K(X^0_{-\infty})=\infty$ if there is no such $K$.

\bigskip
\noindent
{\bf Definition~$2$.} The stationary time series $\{X_n\}$ is said to be finitarily Markovian if 
$K(X^0_{-\infty})$ is finite (though not necessarily bounded) almost surely.

\bigskip
\noindent
In order to estimate $K({\tilde X}^0_{-\infty})$ we need to define some explicit statistics.

\smallskip
\noindent
Define
\begin{eqnarray*}
\lefteqn{
\Delta_k({\tilde X}^0_{-k+1})=}\\
&& \sup_{1\le i} 
\sup_{\{ z^{-k}_{-k-i+1}\in {\cal X}^i, x\in {\cal X}  :   
p_{k+i} (z^{-k}_{-k-i+1},{\tilde X}^0_{-k+1},x)>0 \} }
\left| p_{k-1}(x| {\tilde X}^0_{-k+1})- p_{k+i-1}(x|(z^{-k}_{-k-i+1} ,{\tilde X}^0_{-k+1}))\right|.
\end{eqnarray*}

\noindent
We will divide the data segment $X_0^n$ into two parts: $X_0^{\lceil{n\over 2}\rceil-1}$ and 
$X_{\lceil {n\over 2} \rceil}^n$. Let ${\cal L}_{n,k}^{(1)}$ denote the set of strings with length $k+1$ 
which appear at all in 
$X_0^{\lceil{n\over 2}\rceil-1}$. That is,
$$
{\cal L}_{n,k}^{(1)}= \{x^0_{-k}\in {\cal X}^{k+1}: 
 \exists k\le t \le \lceil{n\over 2}\rceil-1 : X^t_{t-k}=x^0_{-k}\}.
$$ 

\noindent
For a fixed  $0<\gamma<1$  let ${\cal L}_{n,k}^{(2)}$ denote the set of strings with length $k+1$  which appear  
more than $n^{1-\gamma}$ times in $X_{\lceil {n\over 2} \rceil}^n$. That is, 
$${\cal L}_{n,k}^{(2)}=\{x^0_{-k}\in {\cal X}^{k+1}: 
\#\{\lceil {n\over 2} \rceil+k\le t\le n: X^t_{t-k}=x^0_{-k}\} > n^{1-\gamma}\}.$$
Let 
$$
{\cal L}_{k}^n={\cal L}_{n,k}^{(1)}\bigcap {\cal L}_{n,k}^{(2)}.
$$
 We define the empirical version of $\Delta_k$ as follows:  
\begin{eqnarray*}
\lefteqn{ {\hat \Delta}^n_k({\tilde X}^0_{-k+1})=\max_{1\le i \le n}
\max_{(z^{-k}_{-k-i+1},{\tilde X}^0_{-k+1},x)\in {\cal L}^n_{k+i} }
1_{\{\zeta_{J(k)}\le \lceil{n\over 2}\rceil-1\} } }\\
&& \left| 
{ \#\{\lceil {n\over 2} \rceil+k\le t\le n: X^t_{t-k}=({\tilde X}^0_{-k+1},x)\}\over 
\#\{\lceil {n\over 2} \rceil+k-1\le t\le n-1: X^t_{t-k+1}={\tilde X}^0_{-k+1}\}} \right.  \\
&-& \left. 
{ \#\{\lceil {n\over 2} \rceil+k+i\le t\le n: X^t_{t-k-i}=(z^{-k}_{-k-i+1},{\tilde X}^0_{-k+1},x)\}\over 
\#\{\lceil {n\over 2} \rceil+k+i-1\le t\le n-1: X^t_{t-k-i+1}=(z^{-k}_{-k-i+1},{\tilde X}^0_{-k+1})\}}
\right|. 
\end{eqnarray*}
Note that the cut off  $1_{\{\zeta_{J(k)}\le \lceil{n\over 2}\rceil-1\} }$ ensures that ${\tilde X}^0_{-k+1}$ is defined from 
$X_0^{\lceil{n\over 2}\rceil-1}$.

\smallskip
\noindent
Observe, that by ergodicity, for any fixed $k$, 
\begin{equation}\label{Deltatozero}
\liminf_{n\to\infty}{\hat \Delta}^n_k\ge \Delta_k \ \ \mbox{almost surely.}
\end{equation}

\bigskip
\noindent
We define an estimate $\chi_n$ for $K({\tilde X}^0_{-\infty})$  from samples 
$X_0^n$ as follows.  
Let $0< \beta <{1-\gamma \over 2}$ be  arbitrary. Set $\chi_0=0$, and for $n\ge 1$ 
let  $\chi_n$ be the smallest 
$0\le k_n< n $ such that 
${\hat \Delta }^n_{k_n}\le n^{-\beta}$. 

\noindent
Observe  that if $\zeta_j\le \lceil{n\over 2}\rceil-1<\zeta_{j+1}$ then $\chi_n\le l_{j+1}$.

\bigskip
\noindent
Here the idea is (cf. the proof of the  Theorem) 
that if $K({\tilde X}^0_{-\infty})<\infty$ then  $\chi_n$ will be equal to $K({\tilde X}^0_{-\infty})$
 eventually and if $K({\tilde X}^0_{-\infty})=\infty$ then  
 $\chi_n\to\infty$.

\bigskip
\noindent
Now we define the sequence of stopping times $\lambda_n$ along which we will be able to estimate.  
Set $\lambda_0=\zeta_0$, and for $n\ge 1$ 
if $\zeta_j\le \lambda_{n-1}<\zeta_{j+1}$ then put 
\begin{equation}\label{deflambda}
\lambda_n=
\min\{t>\lambda_{n-1} : X_{t-\chi_t+1}^{t}=
X_{\zeta_{j}-\chi_t+1}^{\zeta_{j}}\}
\end{equation}
and
\begin{equation}\label{defkappa}
\kappa_n=\chi_{\lambda_{n}}. 
\end{equation}

\noindent
Observe that if $\zeta_j\le \lambda_{n-1}<\zeta_{j+1}$ then 
$\zeta_j\le \lambda_{n-1}< \lambda_n\le \zeta_{j+1}$. 
If $\chi_{\lambda_{n-1}+1}=0$ then $\lambda_n=\lambda_{n-1}+1$. 
Note that 
$\lambda_n$ is a stopping time and  
$\kappa_n$ is our estimate for $K({\tilde X}^0_{-\infty})$ from samples $X_0^{\lambda_n}$.   

\bigskip
\noindent
Let ${\cal X}^{*-}$ be the set of all one-sided   sequences, that is, 
$${\cal X}^{*-} =\{ (\dots,x_{-1},x_0): x_i\in {\cal X} \ \ \mbox{for all $-\infty<i\le 0$}\}.$$
Let $f: {\cal X}\rightarrow (-\infty,\infty)$ be bounded, 
otherwise arbitrary.  Define the function 
$F : {\cal X}^{*-}\rightarrow (-\infty,\infty)$ 
as 
$$
F(x^{0}_{-\infty})=E(f(X_1)|X^{0}_{-\infty}=x^0_{-\infty}).
$$
E.g. if $f(x)=1_{\{x=z\}}$ for a fixed 
$z\in {\cal X}$ then $F(y^0_{-\infty})=P(X_1=z|X_{-\infty}^0=y^0_{-\infty}).$  If ${\cal X}$ is a finite or 
countably infinite subset of the reals and $f(x)=x$ then 
$F(y^0_{-\infty})=E(X_1|X_{-\infty}^0=y^0_{-\infty}).$

\bigskip
\noindent
One denotes the $n$th estimate of $E(f(X_{\lambda_n+1})|X_0^{\lambda_n})$ 
from samples $X_0^{\lambda_n}$ by $f_n$, 
and defines it to be
\begin{equation}
\label{fkdistrestimate2}
f_n=
{1\over n}\sum_{j=0}^{n-1} f(X_{\lambda_j+1}). 
\end{equation}

\section{Main Results}

\noindent
Define the distance $d^*(\cdot,\cdot)$ on ${\cal X}^{*-}$ as follows. 
For $x^0_{-\infty}$, $y^0_{-\infty}\in {\cal X}^{*-}$ let  
\begin{equation}
\label{defdistance}
d^*(x^0_{-\infty}, y^0_{-\infty})=
\sum_{i=0}^{\infty} 2^{-i-1} 1_{ \{x_{-i} \neq y_{-i}\} }.
\end{equation}

\bigskip
\noindent
{\bf Definition~$3$.}
We say that  $F(X^{0}_{-\infty})$
is almost surely continuous if for some set $C\subseteq {\cal X}^{*-}$  
which has probability one   
the function $F(X^{0}_{-\infty})$ restricted to this set $C$ 
is continuous with respect to metric $ d^*(\cdot,\cdot)$. 
(Cf. Morvai and Weiss \cite{MW03}.)

\bigskip
\noindent
The processes with almost surely continuous conditional expectation 
generalizes the processes for which it is actually continuous, cf.  
Kalikow \cite{Ka90} and Keane \cite{Ke72}. 
The stationary finitarily Markovian processes are included in the class of stationary 
processes with almost surely continuous $E(f(X_1)|X^0_{-\infty})$ for arbitrary bounded $f(\cdot)$.

\bigskip
\noindent
Note that  Ryabko \cite{Ryabko88}, and 
Gy\"orfi, Morvai, Yakowitz \cite{GYMY98} showed that one cannot estimate
$P(X_{n+1}=1|X_{0}^n)$ for all  $n$  in a pointwise consistent way 
even for the class of all stationary and ergodic binary finitarily Markovian 
time series.

\bigskip
\noindent
The entropy rate $H$ associated with a stationary finite or countably infinite 
 alphabet time series
$\{X_n\}$ is defined as
$
H=\lim_{n\to\infty} {-1\over n+1}  \sum_{x_{-n}^0\in {\cal X}^{n+1} }
p_n(x_{-n}^0) \log_2 p_n(x_{-n}^0) $. 
We note that the entropy rate of a stationary finite alphabet time series is finite. 
For details cf. Cover, Thomas \cite{CT91}, pp. 63-64.

\bigskip
\noindent
Fix positive real numbers $0<\beta,\gamma<1$ such that $2\beta+\gamma<1$, fix  a sequence $l_n$ that   
 $1=l_1\le l_2,\dots$, $l_n\to\infty$ and fix  a bounded function 
 $f(\cdot) : {\cal X}\rightarrow (-\infty,\infty)$
 and with these numbers, sequence and function define 
 $\zeta_n$, $\chi_n$, $\kappa_n$, $\lambda_n$ and $F(\cdot)$ as described in the previous section. 
For the resulting $f_n$ we have the following theorem:  

\bigskip
\noindent
{\bf  THEOREM. }
{\it  
Let  $\{X_n\}$ be a stationary and ergodic  time series taking values from a
finite or countably infinite  set ${\cal X}$. 
If the conditional expectation $F(X_{-\infty}^{0})$ is almost surely 
continuous then almost surely, 
$$
\lim_{n\to\infty} f_n=
F({\tilde X}^0_{-\infty}) \ \ \mbox{ and } \ \ 
\lim_{n\to\infty} \left| f_n- E(f(X_{\lambda_n+1})|X_0^{\lambda_n})\right| =0. 
$$

\noindent
The $l_n$ may be chosen in such a fashion that whenever the 
stationary and ergodic time series $\{X_n\}$ has finite entropy rate 
then 
 the 
$\lambda_n$ grow no faster than a polynomial in $n$.

\noindent
If  the stationary and ergodic time series $\{X_n\}$ turns out to be finitarily Markovian then 
$$
\lim_{n\to\infty} {\lambda_n\over n}=
{1\over p_{K({\tilde X}^0_{-\infty})-1}({\tilde X}^0_{-K({\tilde X}^0_{-\infty})+1})}<\infty \ \ \mbox{almost surely}. 
$$ 

\noindent
Moreover, if the stationary and ergodic time series $\{X_n\}$ turns out to be  
independent and identically distributed  then 
$
\lambda_n=\lambda_{n-1}+1 $ 
eventually almost surely.  
}

\bigskip
\noindent
{\bf Proof of the Theorem :}

\noindent
{\bf Step 1.}
{\it The  time series 
$\{{\tilde X}_n\}_{n=-\infty}^0$ and  
$\{X_n\}_{n=-\infty}^{0}$ have 
identical  distribution.  } 

\smallskip
\noindent
For all $k\ge 1$ and 
$1\le i\le k$ 
define (similarly to   Morvai and Weiss \cite{MW03})
$\hat\zeta^k_0=0$ and  
$$
\hat\zeta^k_i=\hat\zeta^k_{i-1}-
\min\{t>0 : 
X_{\hat\zeta^k_{i-1}-(l_{k-i+1}-1)-t}^{\hat\zeta^k_{i-1}-t}
=
{X}_{\hat\zeta^k_{i-1}-(l_{k-i+1}-1)}^{\hat\zeta^k_{i-1}}\}.
$$
Let $T$ denote the left shift operator,  
that is, $(T x^{\infty}_{-\infty})_i=x_{i+1}$. It is easy to see that if 
$\zeta_k(x_{-\infty}^{\infty})=l$ then 
${\hat \zeta}^k_k(T^l x_{-\infty}^{\infty})=-l$.

\noindent
Now the statement follows from stationarity and the fact that 
for $k\ge 0$, $n\ge 0$,   
$x^{0}_{-n}\in {\cal X}^{n+1}$,  $l\ge 0$,  
\begin{equation}\label{shiftequation} 
T^{l} \{X^{\zeta_k}_{\zeta_k-n}=x^0_{-n},\zeta_k=l\} =
\{ X^{0}_{-n}=x^{0}_{-n},{\hat \zeta}^k_k(X^0_{-\infty})=-l\}. 
\end{equation}

\noindent
{\bf Step 2.} {\it We show that 
$P(\chi_n=K({\tilde X}^0_{-\infty}) \  \mbox{eventually}\ |K({\tilde X}^0_{-\infty})<\infty)=1$ and 
$P(\lim_{n\to\infty} \chi_n=\infty |K({\tilde X}^0_{-\infty})=\infty)=1$. 
}

\smallskip
\noindent
By Step~1, $\{{\tilde X}_n\}^0_{n=-\infty}$ is  stationary  and ergodic with the same distribution as $\{X_n\}^0_{n=-\infty}$.
We may assume that the sample path ${\tilde X}^0_{-\infty}$ is such that all finite blocks that appear  have positive probability.
It is immediate that if $K({\tilde X}^0_{-\infty})<\infty$ then 
for all $k\ge K({\tilde X}^0_{-\infty})$, $\Delta_k=0$ and $\Delta_{K({\tilde X}^0_{-\infty})-1}>0$ 
(otherwise the length of the memory  would be not greater than $K({\tilde X}^0_{-\infty})-1$).  
If $K({\tilde X}^0_{-\infty})=\infty$ then $\Delta_k>0$ 
for all $k$, (otherwise $K({\tilde X}^0_{-\infty})$ would be finite).
Thus by (\ref{Deltatozero}) if $K({\tilde X}^0_{-\infty})=\infty$ then $\chi_n\to\infty$ and if 
 $K({\tilde X}^0_{-\infty})<\infty$ then 
$\chi_n\ge K({\tilde X}^0_{-\infty})$ eventually almost surely. 
We have to show that $\chi_n\le K({\tilde X}^0_{-\infty})$ eventually almost surely 
provided that $K({\tilde X}^0_{-\infty})<\infty$. 
 
 \noindent
 Fix now $k<n$. 
We will estimate  the probability of the undesirable  event as follows:
\begin{eqnarray*}
\lefteqn{
P({\hat \Delta}^n_k> n^{-\beta}, K({\tilde X}^0_{-\infty})=k|X_0^{\lceil {n\over 2} \rceil})}\\
&\le&\sum_{i=1}^{n}
  P(\max_{(z^{-k}_{-k-i+1},{\tilde X}^0_{-k+1},x)\in {\cal L}^n_{k+i} }
1_{\{\zeta_{J(k)}\le \lceil{n\over 2}\rceil-1\} } \\
&& \left| 
{ \#\{\lceil {n\over 2} \rceil+k\le t\le n: X^t_{t-k}=({\tilde X}^0_{-k+1},x)\}\over 
\#\{\lceil {n\over 2} \rceil+k-1\le t\le n-1: X^t_{t-k+1}={\tilde X}^0_{-k+1}\}} \right.  \\
&& -  \left. 
{ \#\{\lceil {n\over 2} \rceil+k+i\le t\le n: X^t_{t-k-i}=(z^{-k}_{-k-i+1},{\tilde X}^0_{-k+1},x)\}\over 
\#\{\lceil {n\over 2} \rceil+k+i-1\le t\le n-1: X^t_{t-k-i+1}=(z^{-k}_{-k-i+1},{\tilde X}^0_{-k+1})\}}
\right|\\
&& > n^{-\beta}, K({\tilde X}^0_{-\infty})=k|X_0^{\lceil {n\over 2} \rceil}).
 \end{eqnarray*}

 \noindent
 Define ${\cal M}_{k-1}$ as the set of all $x^0_{-k+1}\in {\cal X}^k$ such that 
 for all $i\ge 1$, $z\in {\cal X}$, and $y^{-k}_{-k-i+1}\in {\cal X}^i$, 
  $p_{k+i}(y^{-k}_{-k-i+1},x^0_{-k+1},z)>0$ implies that 
 $ p_{k-1}(z|x^0_{-k+1})=p_{k+i-1}(z|y^{-k}_{-k-i+1},x^0_{-k+1})
$.
 By the definition of ${\hat \Delta}^n_k$ and since $K({\tilde X}^0_{-\infty})=k$ we have easily that 
 \begin{eqnarray*}
\lefteqn{ P(\max_{(z^{-k}_{-k-i+1},{\tilde X}^0_{-k+1},x)\in {\cal L}^n_{k+i} }
1_{\{\zeta_{J(k)}\le \lceil{n\over 2}\rceil-1\} }} \\
&& \left| 
{ \#\{\lceil {n\over 2} \rceil+k\le t\le n: X^t_{t-k}=({\tilde X}^0_{-k+1},x)\}\over 
\#\{\lceil {n\over 2} \rceil+k-1\le t\le n-1: X^t_{t-k+1}={\tilde X}^0_{-k+1}\}} \right.  \\
&& - \left. 
{ \#\{\lceil {n\over 2} \rceil+k+i\le t\le n: X^t_{t-k-i}=(z^{-k}_{-k-i+1},{\tilde X}^0_{-k+1},x)\}\over 
\#\{\lceil {n\over 2} \rceil+k+i-1\le t\le n-1: X^t_{t-k-i+1}=(z^{-k}_{-k-i+1},{\tilde X}^0_{-k+1})\}}
\right|\\
&& > n^{-\beta}, K({\tilde X}^0_{-\infty})=k|X_0^{\lceil {n\over 2} \rceil})\\
&\le&  
P(\max_{y^0_{-k+1}\in {\cal M}_{k-1},(z^{-k}_{-k-i+1},y^0_{-k+1},x)\in {\cal L}^n_{k+i} } \\
&& \left| 
{ \#\{\lceil {n\over 2} \rceil+k\le t\le n: X^t_{t-k}=(y^0_{-k+1},x)\}\over 
\#\{\lceil {n\over 2} \rceil+k-1\le t\le n-1: X^t_{t-k+1}=y^0_{-k+1}\}} \right.  \\
&& - \left. 
{ \#\{\lceil {n\over 2} \rceil+k+i\le t\le n: X^t_{t-k-i}=(z^{-k}_{-k-i+1},y^0_{-k+1},x)\}\over 
\#\{\lceil {n\over 2} \rceil+k+i-1\le t\le n-1: X^t_{t-k-i+1}=(z^{-k}_{-k-i+1},y^0_{-k+1})\}}
\right|\\
&& > n^{-\beta} |X_0^{\lceil {n\over 2} \rceil}).
\end{eqnarray*}
We can estimate this last probability as the sum of two terms:

\begin{eqnarray*}
\lefteqn{ P(\max_{y^0_{-k+1}\in {\cal M}_{k-1},(z^{-k}_{-k-i+1},y^0_{-k+1},x)\in {\cal L}^n_{k+i} }
} \\
&& \left| 
{ \#\{\lceil {n\over 2} \rceil+k\le t\le n: X^t_{t-k}=(y^0_{-k+1},x)\}\over 
\#\{\lceil {n\over 2} \rceil+k-1\le t\le n-1: X^t_{t-k+1}=y^0_{-k+1}\}} \right.  \\
&& - \left. 
{ \#\{\lceil {n\over 2} \rceil+k+i\le t\le n: X^t_{t-k-i}=(z^{-k}_{-k-i+1},y^0_{-k+1},x)\}\over 
\#\{\lceil {n\over 2} \rceil+k+i-1\le t\le n-1: X^t_{t-k-i+1}=(z^{-k}_{-k-i+1},y^0_{-k+1})\}}
\right| > n^{-\beta} |X_0^{\lceil {n\over 2} \rceil})\\
&\le& 
P(\max_{y^0_{-k+1}\in {\cal M}_{k-1},(z^{-k}_{-k-i+1},y^0_{-k+1},x)\in {\cal L}^n_{k+i} }
 \\
&& \left| 
{ \#\{\lceil {n\over 2} \rceil+k\le t\le n: X^t_{t-k}=(y^0_{-k+1},x)\}\over 
\#\{\lceil {n\over 2} \rceil+k-1\le t\le n-1: X^t_{t-k+1}=y^0_{-k+1}\}} - 
p_{k-1}(x|y^0_{-k+1})
\right|\\
 && >0.5 n^{-\beta}|X_0^{\lceil {n\over 2} \rceil})\\
&+& 
P(\max_{y^0_{-k+1}\in {\cal M}_{k-1},(z^{-k}_{-k-i+1},y^0_{-k+1},x)\in {\cal L}^n_{k+i} }
 \\
&& \left| 
p_{k-1}(x|y^0_{-k+1})-
{ \#\{\lceil {n\over 2} \rceil+k+i\le t\le n: X^t_{t-k-i}=(z^{-k}_{-k-i+1},y^0_{-k+1},x)\}\over 
\#\{\lceil {n\over 2} \rceil+k+i-1\le t\le n-1: X^t_{t-k-i+1}=(z^{-k}_{-k-i+1},y^0_{-k+1})\}}\right|\\
&& >0.5 n^{-\beta}|X_0^{\lceil {n\over 2} \rceil}).
\end{eqnarray*}
We overestimate these probabilities. 
For any $m\ge 0$ and $x^0_{-m}$ define  $\sigma^m_i(x^0_{-m})$ as the time of the $i$-th ocurrence of the string $x^0_{-m}$ in 
 the data segment $X_{\lceil {n\over 2} \rceil}^n$, that is, let $\sigma^m_{0}(x^0_{-m})=\lceil {n\over 2} \rceil+m-1$ 
 and for $i\ge 1$ define
 $$
 \sigma^m_i(x^0_{-m})=\min\{t>\sigma^m_{i-1}(x^0_{-m}) : X^t_{t-m}=x^0_{-m}\}.
 $$
 Now  
 \begin{eqnarray*}
\lefteqn{ P(\max_{y^0_{-k+1}\in {\cal M}_{k-1},(z^{-k}_{-k-i+1},y^0_{-k+1},x)\in {\cal L}^n_{k+i} }
} \\
&& \left| 
{ \#\{\lceil {n\over 2} \rceil+k\le t\le n: X^t_{t-k}=(y^0_{-k+1},x)\}\over 
\#\{\lceil {n\over 2} \rceil+k-1\le t\le n-1: X^t_{t-k+1}=y^0_{-k+1}\}} \right.  \\
&& - \left. 
{ \#\{\lceil {n\over 2} \rceil+k+i\le t\le n: X^t_{t-k-i}=(z^{-k}_{-k-i+1},y^0_{-k+1},x)\}\over 
\#\{\lceil {n\over 2} \rceil+k+i-1\le t\le n-1: X^t_{t-k-i+1}=(z^{-k}_{-k-i+1},y^0_{-k+1})\}}
\right| > n^{-\beta} |X_0^{\lceil {n\over 2} \rceil})\\
&\le& 
 P(\max_{y^0_{-k+1}\in {\cal M}_{k-1},(y^0_{-k+1},x)\in {\cal L}^{(1)}_{n,k} }\sup_{j> n^{1-\gamma}}
 \\
&& \left| {1\over j} \sum_{r=1}^j 1_{\{X_{\sigma_r^{k-1}(y^0_{-k+1})}=x\}} -p_{k-1}(x|y^0_{-k+1})\right| 
  >0.5 n^{-\beta}|X_0^{\lceil {n\over 2} \rceil})\\
&+&
  P(\max_{y^0_{-k+1}\in {\cal M}_{k-1},(z^{-k}_{-k-i+1},y^0_{-k+1},x)\in {\cal L}^{(1)}_{n,k+i} }\sup_{j> n^{1-\gamma}} \\
&& \left| {1\over j} \sum_{r=1}^j 1_{\{X_{\sigma_r^{k+i-1}(z^{-k}_{-k-i+1},y^0_{-k+1},x)}=x\}} -
p_{k-1}(x|y^0_{-k+1})\right| 
  >0.5 n^{-\beta}|X_0^{\lceil {n\over 2} \rceil})
\end{eqnarray*}
Since both ${\cal L}^{(1)}_{n,k}$ and ${\cal L}^{(1)}_{n,k+i}$ depend solely on $X_0^{\lceil {n\over 2} \rceil}$ we get 
\begin{eqnarray*}
\lefteqn{ P(\max_{y^0_{-k+1}\in {\cal M}_{k-1},(z^{-k}_{-k-i+1},y^0_{-k+1},x)\in {\cal L}^n_{k+i} }
} \\
&& \left| 
{ \#\{\lceil {n\over 2} \rceil+k\le t\le n: X^t_{t-k}=(y^0_{-k+1},x)\}\over 
\#\{\lceil {n\over 2} \rceil+k-1\le t\le n-1: X^t_{t-k+1}=y^0_{-k+1}\}} \right.  \\
&& - \left. 
{ \#\{\lceil {n\over 2} \rceil+k+i\le t\le n: X^t_{t-k-i}=(z^{-k}_{-k-i+1},y^0_{-k+1},x)\}\over 
\#\{\lceil {n\over 2} \rceil+k+i-1\le t\le n-1: X^t_{t-k-i+1}=(z^{-k}_{-k-i+1},y^0_{-k+1})\}}
\right| > n^{-\beta} |X_0^{\lceil {n\over 2} \rceil})\\
&\le& \sum_{y^0_{-k+1}\in {\cal M}_{k-1},(y^0_{-k+1},x)\in {\cal L}^{(1)}_{n,k} }
\sum_{j= \lceil n^{1-\gamma}\rceil}^{\infty}
   P(
\left| {1\over j} \sum_{r=1}^j 1_{\{X_{\sigma_r^{k-1}(y^0_{-k+1})}=x\}} -p_{k-1}(x|y^0_{-k+1})\right|\\
&&  >0.5 n^{-\beta}|X_0^{\lceil {n\over 2} \rceil})\\
  &+&\sum_{y^0_{-k+1}\in {\cal M}_{k-1},(z^{-k}_{-k-i+1},y^0_{-k+1},x)\in {\cal L}^{(1)}_{n,k+i} }
  \sum_{j= \lceil n^{1-\gamma}\rceil}^{\infty}
   P( \left| {1\over j} \sum_{r=1}^j 1_{\{X_{\sigma_r^{k+i-1}(z^{-k}_{-k-i+1},y^0_{-k+1})}=x\}} \right.\\
  &&  \left. -p_{k-1}(x|y^0_{-k+1})
   \right|
  >0.5 n^{-\beta}|X_0^{\lceil {n\over 2} \rceil}).
  \end{eqnarray*}
   Each of these represents the deviation of an empirical count from its mean. 
The variables in question are independent since whenever the block 
 $y^0_{-k+1}$ occurs the next term is chosen using the same distribution $p_{k-1}(x|y^0_{-k+1})$.
 Thus by Hoeffding's inequality (cf. Hoeffding \cite{Hoeffding63} or Theorem 8.1 of Devroy et. al. \cite{DGYL96}) for 
 sums of bounded independent random variables and since the cardinality of both 
 ${\cal L}^{(1)}_{n,k}$ and ${\cal L}^{(1)}_{n,k+i}$ is not greater than $ (n+2)/2$, 
we have 
\begin{eqnarray*}
\lefteqn{ P(\max_{y^0_{-k+1}\in {\cal M}_{k-1},(z^{-k}_{-k-i+1},y^0_{-k+1},x)\in {\cal L}^n_{k+i} }
} \\
&& \left| 
{ \#\{\lceil {n\over 2} \rceil+k\le t\le n: X^t_{t-k}=(y^0_{-k+1},x)\}\over 
\#\{\lceil {n\over 2} \rceil+k-1\le t\le n-1: X^t_{t-k+1}=y^0_{-k+1}\}} \right.  \\
&-& \left. 
{ \#\{\lceil {n\over 2} \rceil+k+i\le t\le n: X^t_{t-k-i}=(z^{-k}_{-k-i+1},y^0_{-k+1},x)\}\over 
\#\{\lceil {n\over 2} \rceil+k+i-1\le t\le n-1: X^t_{t-k-i+1}=(z^{-k}_{-k-i+1},y^0_{-k+1})\}}
\right| > n^{-\beta} |X_0^{\lceil {n\over 2} \rceil})\\
&\le & 2 {n+2\over 2} \sum_{j= \lceil n^{1-\gamma}\rceil}^{\infty} 2 e^{-2n^{-2\beta}j}.
\end{eqnarray*}
Thus 
\begin{eqnarray*}
\lefteqn{P({\hat \Delta}^n_k> n^{-\beta},K({\tilde X}^0_{-\infty})=k|X_0^{\lceil {n\over 2} \rceil})}\\
&\le&  n (n+2)   2 e^{-2n^{-2\beta+1-\gamma}}.
\end{eqnarray*}
Integrating both sides we get 
\begin{eqnarray*}
\lefteqn{
P({\hat \Delta}^n_k> n^{-\beta}, K({\tilde X}^0_{-\infty})=k)}\\
&\le&  n (n+2)   2 e^{-2n^{-2\beta+1-\gamma}}.
\end{eqnarray*}
The right hand side is summable provided $2\beta+\gamma<1$ and 
the Borel-Cantelli Lemma yields that  
\begin{eqnarray*}
\lefteqn{
P({\hat \Delta}^n_k\le n^{-\beta}  eventually,K({\tilde X}^0_{-\infty})=k)}\\
&=&P(K({\tilde X}^0_{-\infty})=k).
\end{eqnarray*}
Thus  
$\chi_n\le  k$ eventually almost surely on $ K({\tilde X}^0_{-\infty})=k$.

\noindent
{\bf Step 3.}
{\it We show the first part of the Theorem. }

\noindent
Recalling~(\ref{fkdistrestimate2}) we can write
\begin{equation}
\label{decomposition} 
f_n = {1\over n}\sum_{j=0}^{n-1} [  f(X_{\lambda_j+1})  -
E(f(X_{\lambda_j+1})|X_{-\infty}^{\lambda_j})] + 
 {1\over n}\sum_{j=0}^{n-1} E(f(X_{\lambda_j+1})|X_{-\infty}^{\lambda_j})
 \end{equation}

\noindent
Observe that the first term is an average of orthogonal bounded random variables and 
 by Theorem~3.2.2 in R\'ev\'esz \cite{Revesz68}, it tends to zero. 

\noindent
Now we deal with the second term. 
If $K({\tilde X}^0_{-\infty})<\infty$  then  by Step~2, 
${\chi}_n= K({\tilde X}^0_{-\infty})$ eventually 
 and by 
(\ref{defzeta}), (\ref{defprocesses}), (\ref{deflambda}) and  
Step~1, eventually,   
$$E(f(X_{\lambda_j+1})|X_{-\infty}^{\lambda_j})=
E(f(X_{\lambda_j+1})|X_0^{\lambda_j})=
F({\tilde X}^0_{-\infty}).$$ 

\noindent
We may deal with the case when $K({\tilde X}^0_{-\infty})=\infty$ 
and by Step~2, ${\chi}_n\to\infty$.   
For  arbitrary $j\ge 0$, by ~(\ref{defkappa}) and ~(\ref{deflambda}) and  
the construction in (\ref{defprocesses}),     
\begin{equation}\label{firstjbitequal}
X_{\lambda_j-\kappa_j+1}^{\lambda_j}=
{\tilde X}^0_{-\kappa_j+1} \ \ \mbox{and} \ \  
\lim_{j\to\infty} 
d^*({\tilde X}^0_{-\infty},X_{-\infty}^{\lambda_j})=0  \ \ \mbox{almost surely.}
\end{equation}
Be Step~1, and the almost sure continuity of $F(\cdot)$, for some set  $C\subseteq {\cal
X}^{*-}$ with full measure,    $F(\cdot)$ is continuous on $C$ and   
\begin{equation}\label{inCtilde}
{\tilde X}^{0}_{-\infty}\in C, X_{-\infty}^n\in C \ \ \mbox{for all $n\ge 0$ almost surely.} 
\end{equation}  
By  the continuity of $F(\cdot)$ on the set $C$ and (\ref{firstjbitequal}),   
$E(f(X_{\lambda_j+1})|X^{\lambda_j}_{-\infty})=F(X_{-\infty}^{\lambda_j})
\to F({\tilde X}^{0}_{-\infty})$
 and 
$f_n\to F({\tilde X}^{0}_{-\infty})$ almost surely.

\noindent 
Define the random neighbourhood ${\cal  N}_j(X_{0}^{\lambda_j})$ of $X_{0}^{\lambda_j}$ 
depending on the random data segment  
$X_{0}^{\lambda_j}$ itself as  
$$
{\cal  N}_j(X_{0}^{\lambda_j})=\{ z^0_{-\infty}\in {\cal X}^{*-} \ : \ 
z_{-{\kappa_j}+1}= X_{\lambda_j-{\kappa_j}+1},\dots, 
z_{0}=X_{\lambda_j} \}.   
$$

\noindent
Note that by (\ref{defzeta}), (\ref{defprocesses}), (\ref{defkappa}) and (\ref{deflambda}),   
$
{\tilde X}^0_{-\infty}\in {\cal  N}_j(X_{0}^{\lambda_j})  
$ and 
by (\ref{inCtilde}) and 
the continuity of $F(\cdot)$ on the set $C$, and since $\kappa_j\to\infty$,   
by~(\ref{firstjbitequal}), 
 almost surely, 
\begin{eqnarray*}
\lefteqn{ \lim_{j\to\infty} 
\left|E(f(X_{\lambda_j+1})|X^{\lambda_j}_{0})-F({\tilde X}^0_{-\infty})\right|
=
\lim_{j\to\infty} 
\left|E\{ F(X^{\lambda_j}_{-\infty}) |X^{\lambda_j}_0\}-F({\tilde X}^0_{-\infty})\right|}\\
&\le&
\lim_{j\to\infty} 
\sup_{y^0_{-\infty}, z^0_{-\infty}\in {\cal  N}_j(X_{0}^{\lambda_j})
\bigcap C}
|F(y^0_{-\infty})-F(z^0_{-\infty})|=0.
\end{eqnarray*}

\noindent
{\bf Step 4.}
{\it We show the second part of the Theorem.}

\noindent
Now we assume that  the stationary and ergodic finite or countably infinite alphabet 
time series $\{X_n\}$ possesses  finite entropy rate $H$. 
(A stationary finite alphabet time series always has finite entropy rate.) 

\noindent
We will in fact obtain a more precise estimate, namely, if 
 for some $0<\epsilon_2<\epsilon_1$,   
$\sum_{k=1}^{\infty} (k+1) 2^{- l_k(\epsilon_1-\epsilon_2)}<\infty$ 
then 
$$
\lambda_n< 2^{l_n(H+\epsilon_1)} 
\ \ \mbox{ eventually almost surely.} 
$$
In particular,  for arbitrary $\delta>0$, $0<\epsilon_2<\epsilon_1$, 
if
$l_n=\min\left(n,\max\left(1, 
\lfloor {2+\delta \over \epsilon_1-\epsilon_2 } \log_2 n\rfloor\right) \right)$ 
then  
$$\lambda_n< n^{ {2+\delta\over \epsilon_1-\epsilon_2} (H+\epsilon_1) }$$
 eventually almost surely, and 
the upper bound is a polynomial.

\bigskip
\noindent
Since $\lambda_n\le \zeta_n$, it is enough to prove the result for $\zeta_n$.  
Let ${\cal X}^*$ be the set of all two-sided   sequences, that is, 
$${\cal X}^*= \{ (\dots,x_{-1},x_0,x_1,\dots):
x_i\in {\cal X} \ \ \mbox{for all $-\infty\le i< \infty$}\}.
$$

\noindent
Define $B_k\subseteq {\cal X}^{l_k}$ as 
$
B_k=\{ x^0_{-l_k+1} \in {\cal X}^{l_k}: 2^{-l_k(H+\epsilon_2)} < p_{l_k-1}(x^0_{-l_k+1})\}.
$
Note that there is a trivial bound on the cardinality of the set $B_k$, namely,
\begin{equation} \label{cardbound}
|B_k|\le 2^{l_k(H+\epsilon_2)}. 
\end{equation}

\noindent
Define the set $\Upsilon_k(y^0_{-k+1})$ as follows: 
$$
\Upsilon_k(y^0_{-l_k+1})=\{z^{\infty}_{-\infty}\in {\cal X}^+ :  
-{\hat \zeta}^k_k(z_{-\infty}^0) \ge 2^{l_k(H+\epsilon_1)}, 
z^0_{-l_k+1}=y^0_{-l_k+1})\}.
$$
We will estimate the probability of $\Upsilon_k(y^0_{-l_k+1})$  by 
a frequency argument. 
Let $x^{\infty}_{-\infty}\in {\cal X}^*$ be a typical sequence of the time series $\{X_n\}$. 
Define $\rho_0(y^0_{-l_k+1},x_{-\infty}^{\infty})=0$ 
and for $i\ge 1$ let 
$$
\rho_i(y^0_{-l_k+1},x_{-\infty}^{\infty})=
\min \{l> \rho_{i-1}(y^0_{-l_k+1},x_{-\infty}^{\infty}): 
T^{-l} x_{-\infty}^{\infty}\in \Upsilon_k(y^0_{-l_k+1})\}.
$$
Define also $\tau_0(y^0_{-l_k+1},x_{-\infty}^{\infty})=0$  
and for $i\ge 1$ let 
$$
\tau_i(y^0_{-l_k+1},x_{-\infty}^{\infty})=\min \{l\ge  
\tau_{i-1}(y^0_{-l_k+1},x_{-\infty}^{\infty})
+2^{l_k(H+\epsilon_1)}: 
T^{-l} x_{-\infty}^{\infty}\in \Upsilon_k(y^0_{-l_k+1})\}.
$$
Notice that if  $\tau_{i-1}=\rho_m$ then $ \tau_i\le \rho_{m+k+1}$. 
(Indeed, since there are at least $k+1$ occurrences of the block $y^0_{-l_k+1}$ 
in the data segment
$X_{-\rho_{m+k+1}-l_k+1}^{\rho_m+1}$ hence 
$2^{l_k(H+\epsilon_1)}\le 
-{\hat \zeta}_k^k(T^{-\rho_m} x_{-\infty}^{\infty})\le \rho_{m+k+1}-\tau_{i-1}$.)
By the ergodicity of the time series $\{X_n\}$, 
\begin{eqnarray}
\lefteqn{ \nonumber
P( X_{-\infty}^{\infty} \in \Upsilon_k(y^0_{-l_k+1}) )=
\lim_{t\to\infty} 
{ \# \{ j\ge 1 :   
\rho_j(y^0_{-l_k+1},x_{-\infty}^{\infty})\le \tau_t(y^0_{-l_k+1},x_{-\infty}^{\infty})\} 
\over \tau_t(y^0_{-l_k+1},x_{-\infty}^{\infty})} } \\
&=&  \nonumber\lim_{t\to \infty} 
{\sum_{l=1}^{t} 
 \# \{ j\ge 1 :  \tau_{l-1}(y^0_{-l_k+1},x_{-\infty}^{\infty}) < 
\rho_j(y^0_{-l_k+1},x_{-\infty}^{\infty}) \le \tau_l(y^0_{-l_k+1},x_{-\infty}^{\infty})\}
\over \tau_t(y^0_{-l_k+1},x_{-\infty}^{\infty})} 
 \\
&\le&  \lim_{t\to\infty} {t (k+1)\over t 2^{l_k(H+\epsilon_1)} }
= \label{upperbound} {(k+1)\over 2^{l_k(H+\epsilon_1)}}. 
\end{eqnarray}
Since 
$$
T^l\{\zeta_k=l,  X^{\zeta_k}_{\zeta_k-l_k+1}\in B_k\}=
\{{\hat \zeta}_k^k=-l, X^0_{-l_k+1}\in B_k\}
$$
by stationarity  and  the upper bound on the cardinality of the set $B_k$   
in (\ref{cardbound}) and by (\ref{upperbound}), 
   we get  
\begin{eqnarray}
\nonumber
 P(\zeta_k\ge 2^{l_k(H+\epsilon_1)},  {\tilde X}^{0}_{-l_k+1}\in B_k)
&=&
P(\zeta_k\ge 2^{l_k(H+\epsilon_1)},  X^{\zeta_k}_{\zeta_k-l_k+1}\in B_k)\\
&=&\nonumber
P(-{\hat \zeta}_k^k\ge 2^{l_k(H+\epsilon_1)}, X^0_{-l_k+1}\in B_k)\\ 
&=& \nonumber
\sum_{y^0_{-l_k+1}\in B_k} 
P(X_{-\infty}^{\infty} \in \Upsilon_k(y^0_{-l_k+1}))\\
&\le& \nonumber
(k+1) 2^{-l_k (\epsilon_1-\epsilon_2)}.
\end{eqnarray}
By assumption,  the right hand side  sums and  
the Borel-Cantelli Lemma yields 
that the event 
$\{ \zeta_k\ge 2^{l_k(H+\epsilon_1)},{\tilde X}^{0}_{-l_k+1}\in B_k\}$ cannot happen infinitely 
many times. 
By Step 1, the distribution of the time series 
$\{\tilde X_n\}$ is the same as the distribution of $\{X_n\}$ 
and  by the Shannon-McMillan-Breiman Theorem   (cf. Chung \cite{Chung61})
$ {\tilde X}^0_{-l_k+1}\in B_k$ eventually almost surely and so 
$\zeta_k\ge 2^{l_k(H+\epsilon_1)}$ cannot happen infinitely 
many times. 

\noindent
{\bf Step 5.} 
{\it We show the rest of the Theorem.}

\noindent
By Step 2,   if $1\le K({\tilde X}^0_{-\infty})<\infty$ then  
 $\chi_n=K({\tilde X}^0_{-\infty})$ eventually, and by ergodicity, 
${n\over \lambda_n}\to p_{K({\tilde X}^0_{-\infty})-1}({\tilde X}_{-K({\tilde X}^0_{-\infty})+1}^0)>0$. 
If $K({\tilde X}^0_{-\infty})=0$ then 
by Step 2, $\chi_n=0$ eventually, and by (\ref{deflambda}), $\lambda_n=\lambda_{n-1}+1$ 
eventually. The proof of the Theorem is complete.


\end{document}